%% file: EuCAP2023_template.tex
\documentclass[conference,a4paper]{IEEEtran}
\usepackage{xcolor}
\usepackage{balance}

\ifCLASSINFOpdf
 \usepackage[pdftex]{graphicx}
\else
 \usepackage[dvips]{graphicx}
\fi
\ifCLASSOPTIONcompsoc
 \usepackage[caption=false,font=normalsize,labelfont=sf,textfont=sf]{subfig}
\else
 \usepackage[caption=false,font=footnotesize]{subfig}
\fi
\newcommand{\td}[1]{\todo[inline,size=\small]{#1}} 
\usepackage{cite}
\usepackage{amsmath,amssymb,amsfonts}
\usepackage{algorithmic}
\usepackage{graphicx}
\usepackage{textcomp}
\usepackage{xcolor}
\usepackage{comment}
\usepackage{todonotes}
\usepackage{mathtools}
\usepackage{cleveref}
\usepackage{siunitx}

\input{Macros.tex}
\input{Packages.tex}

\begin{document}
%
\title{Fast Direct Solvers for Integral Equations at Low-Frequency Based on Operator Filtering}

\author{\IEEEauthorblockN{
Clément~Henry\IEEEauthorrefmark{2},
Davide~Consoli\IEEEauthorrefmark{1},
Alexandre~Dély\IEEEauthorrefmark{1},
Lyes~Rahmouni\IEEEauthorrefmark{1},
Adrien~Merlini\IEEEauthorrefmark{2},
Francesco~P.~Andriulli\IEEEauthorrefmark{1}
}   
\IEEEauthorblockA{\IEEEauthorrefmark{1}
Department of Electronics and Telecommunications, Politecnico di Torino, Turin, Italy, francesco.andriulli@polito.it}
\IEEEauthorblockA{\IEEEauthorrefmark{2}
Microwave Department, IMT Atlantique, Brest, France, clement.henry@imt-atlantique.fr}
}




\maketitle

\begin{abstract}
This paper focuses on fast direct solvers for integral equations in the low-to-moderate-frequency regime obtained by leveraging preconditioned first kind or second kind operators regularized with Laplacian filters. 
The spectral errors arising from boundary element discretizations are properly handled by filtering that, in addition, allows for the use of low-rank representations for the compact perturbations of all operators involved. Numerical results show the effectiveness of the approaches and their effectiveness in the direct solution of integral equations.
\end{abstract}

\vskip0.5\baselineskip

%

\section{Introduction}

Several inverse strategies require the solutions of multiple forward problems that, with subsequent optimizations, converge to the desired result. In this context, direct solutions (as opposed to iterative ones) are preferable because they are favorably performing when dealing with multiple right-hand-side problems.

Among forward solution strategies, the boundary element method (BEM) is a powerful technique to model time-harmonic electromagnetic scattering from arbitrary geometries. Its advantages are, among others, that only the scatterer’s boundary needs to be discretized and that it automatically enforces the radiation conditions. A popular choice as a forward solution equation is the Electric Field Integral Equation (EFIE). Although this formulation is widespread, its numerical solution comes with its own set of challenges: the linear systems that stem from its discretization are dense and ill-conditioned, with a condition number that grows with the inverses of the frequency and the average mesh edge length \cite{adrianElectromagneticIntegralEquations2021}.

On the one hand, fast algorithms such as the fast multipole method are commonly used to obtain a solution in linear or quasi-linear time and memory complexities. On the other hand, the conditioning problem is usually tackled by Calderón preconditioning \cite{andriulli_well-conditioned_2013}, multilevel schemes, algebraic approaches, etc. In fact, the classical approach is to separate the problem of acceleration from the one of regularization in order to fix them independently.

This work will present a different paradigm, proposing a new family of techniques that combine these two aspects by taking advantage of the common structure present in fast solvers and preconditioners for the matrices of discretized boundary integral operators. In particular, by an operation of regularization we obtain a formulation from which a fast solver naturally emerges. This effect, however, cannot be obtained in a straightforward manner because of the spectral pollution of standard integral operators. A key enabling factor will thus be a recently introduced tool: the quasi-Helmholtz  Laplacian Filter that will be able to fix the spectral behavior of the operators involved allowing for their fast direct inversion with a focus on problems from low to moderate frequencies.

Numerical results will corroborate the theory, showing the practical relevance of all newly proposed techniques.




\section{Notation and Background}

Consider a smooth, closed curve $\Gamma_\mathrm{2}\subset\R^2$ modeling a PEC scatterer in a medium of permeability $\mu$, permittivity $\epsilon$, wavenumber $k=\omega \sqrt{\mu \epsilon}$, and impedance $\eta = \sqrt{\mu/\epsilon}$. The TE electric and magnetic field integral equations (TE-EFIE and TE-MFIE) that relate the tangential component of the current $j_t$ with the tangential component of the incident electric field $e_t^\mathrm{inc}$ and the out-of-plane component of the incident magnetic field $h_z^\mathrm{inc}$, respectively, are
\begin{align} 
    (\im k)^{-1} \op N j_t &= - \eta^{-1} e_t^\mathrm{inc}  \label{eq:TE-EFIE}\\
    \left(\op I/2 -  \op D \right) j_t &= - h_z^\mathrm{inc}  \label{eq:TE-MFIE}
\end{align}
with 
\begin{align}
\left( \op N j_t \right)(\vr) &\coloneqq -\frac{\partial}{\partial n} \int_{\Gamma_2} \frac{\partial}{\partial n'} g_\mathrm{2}(\vr, \vr') j_t(\vr') \diff \vr'\\
\left( \op D j_t \right)(\vr) &\coloneqq  \int_{\Gamma_2} \frac{\partial}{\partial n'} g_\mathrm{2}(\vr, \vr') j_t(\vr') \diff \vr'\, ,
\end{align}
and $g_\mathrm{2}(\vr,\vr') \coloneqq \im/4 H_0^{1}(k|\vr-\vr'|)$. Leveraging the Calderón formula $\op S \op N = \op I / 4 + \op C_2$, where $\op I$ is the identity, $\op C_2$ is a compact operator, and $\left(\op S j_t \right)(\vr) \coloneqq  \int_{\Gamma_2}  g_\mathrm{2}(\vr, \vr') j_t(\vr') \diff \vr'$, \eqref{eq:TE-EFIE} can be turned into a second kind integral equation. 

After expanding the unknown current density with  piecewise linear Lagrange interpolants $\{\varphi_i\}$, defined on a mesh of $\Gamma_2$ made out of segments of average length $h$, as $j_t \approx \sum_{i=1}^N [\vec j_2]_i \varphi_i$ and after Galerkin testing, the discretized forms of the normalized Calderón TE-EFIE and normalized TE-MFIE are
\begin{gather} 
    \mat Z_{2} \vec j_2 = \vec v_{e,2} \label{eq: DiscCaldTEEFIE}\,, \\
    (\mat I /2 - \mathh^{-\frac{1}{2}} \mat D \mathh^{-\frac{1}{2}}) \vec j_2 = \vec v_{h,2} \,,\label{eq: DiscTEEFIE}
\end{gather}
where 
$\mat Z_{2} \coloneqq (\im k)^{-1} \mathh^{-\frac{1}{2}} \mat S \mathh^{-1} \mat N \mathh^{-\frac{1}{2}}$,
$\vec v_{e,2} \coloneqq - \eta^{-1} \mathh^{-\frac{1}{2}} \mat S \mathh^{-1} \vec e_2$, $\vec v_{h,2} \coloneqq - \mathh^{-\frac{1}{2}} \vec h_2$,
$[\mat N]_{ij}\coloneqq\langle \varphi_i , \op N(\varphi_j) \rangle$, $[\mat S]_{ij} \coloneqq \langle \varphi_i , \op S(\varphi_j) \rangle$, $[\mat D]_{ij} \coloneqq \langle \varphi_i , \op D(\varphi_j) \rangle$, $[\vec e_2] _i \coloneqq \langle \varphi_i, e_t^\mathrm{inc} \rangle$, $[\vec h_2] _i \coloneqq \langle \varphi_i, h_z^\mathrm{inc} \rangle$, $[\mathh]_{ij} \coloneqq \langle \varphi_i , \varphi_j \rangle$, and $\langle a, b\rangle \coloneqq \int_{\Gamma_2} a b \diff s$.

In the following, we also consider the 3-dimensional counterpart of the problem described above. Consider a closed, simply-connected domain $\Omega \subset \R^3$ delimited by smooth boundary $\Gamma_3$ and associated outgoing unit surface normal vector field $\uv{n}$. The electric current density $\vt j$ induced on the object by impinging electric and magnetic fields $\vt{e}^\mathrm{inc}$ and $\vt{h}^\mathrm{inc}$ is obtained by solving the 3D EFIE or MFIE
\begin{gather} 
 \vecop T \vt j \coloneqq \im k \vecop T_\mathrm{A} \vt j + (\im k)^{-1} \vecop T_\mathrm{\Phi} \vt j  = - \uv{n} \times\vt e^\mathrm{inc}\,, \label{eq:EFIE} \\
 \vecop M \vt j \coloneqq \left( \vecop I /2 + \vecop K \right) \vt j  = \uv{n} \times\vt h^\mathrm{inc}\,, \label{eq:MFIE}
\end{gather}	 
where $\left( \vecop T_\mathrm{A} \vt j \right) (\vr):= \uv{n}(\vr) \times \int_{\Gamma_3} g_3(\vr,\vr') \vt j(\vr') \diff \vr'$, $\left( \vecop T_\mathrm{\Phi} \vt j \right)(\vr):= - \uv{n}(\vr) \times \nabla_{\vr}\int_{\Gamma_3} g_3(\vr,\vr') \nabla_{\vr'} \cdot \vt j(\vr') \diff \vr' $, $\left( \vecop K \vt j \right) (\vr):= - \uv{n}(\vr) \times \int_{\Gamma_3} \nabla_r g_3(\vr,\vr') \times \vt j(\vr') \diff \vr'$, and $g_3(\vr,\vr') \coloneqq \exp{(\im k |\vr - \vr'|)}/|\vr - \vr'|$ is the free-space Green's function. Leveraging the identity $\vecop T^2 = -\vecop I  / 4 + \vecop C_3 $, a second kind integral equation can be obtained from the EFIE, yielding the standard Calderón EFIE \cite{andriulli_well-conditioned_2013}
\begin{equation} \label{eq:CaldEFIE}
 -\vecop T^2 \vt j = \vecop T \left( \uv{n} \times\vt{e}^\mathrm{inc} \right)\,.
\end{equation}	
To discretize \eqref{eq:EFIE} following a Petrov-Galerkin approach, we first expand the unknown current density as $\vt j \approx \sum_{i=1}^N [\vec j_3]_i \vt{f}_i$, where $\{\vt{f}_i\}$ is the set of Rao-Wilton-Glisson (RWG) basis functions defined on a triangular discretization of $\Gamma_3$ in which the average edge length of the triangles is denoted $h$. We also define on this discretized geometry the sets of rotated RWG functions $\{ \uv{n} \times \vt{f}_i \}$, Buffa-Christiansen (BC) functions $\{\vt{\tilde{f}}_i\}$ \cite{andriulli_well-conditioned_2013}, rotated BC functions $\{\uv{n} \times\vt{\tilde{f}}_i\}$, pyramid functions $\{\lambda_i\}$, and patch functions $\{p_i\}$. We also define the (potentially mixed) Gram matrices for these bases as
$[\mat G_{\vt u \vt v}]_{ij} \coloneqq \langle \vt u_i, \vt v_j\rangle$, and $[\mat {\tilde{G}}_{\vt u \vt v}]_{ij} \coloneqq \langle \uv{n} \times \vt{u}_i, \vt{v}_j\rangle$.
Finally, upon adequate Petrov-Galerkin testing, the discrete normalized Calderón EFIE \eqref{eq:CaldEFIE} and normalized MFIE \eqref{eq:MFIE} are
\begin{gather} 
   \mat Z_3 \vec j_3 = \vec v_{e,3}  \label{eq:DiscCaldEFIE}\,,\\
   \left( \mat I /2 + \mat{\tilde{K}} \right) \vec j_3 = \vec v_{h,3} \,,
\end{gather}
where
\begin{gather}
    \mat Z_3 = -\matRR^{\frac{1}{2}} \matMixTr^{-1} \mat{T_{BC}} \matMix^{-1} \mat{T} \matRR^{-\frac{1}{2}}\,, \\
    \mat{\tilde{K}} = \matRR^{\frac{1}{2}} \matMixTr^{-1} \mat K  \matRR^{-\frac{1}{2}} \,, \\
    \vec v_{e,3} = \matRR^{\frac{1}{2}} \matMixTr^{-1} \mat{T_{BC}} \matMix^{-1} \vec e_3\,,\\
    \vec v_{h,3} = \matRR^{\frac{1}{2}} \matMixTr^{-1}  \vec h_3\,,
\end{gather}
$[\mat T]_{ij}\coloneqq\langle \uv n \times \vt{f}_i , \vecop T(\vt{f}_j) \rangle$, $[\mat K]_{ij}\coloneqq\langle \uv n \times \vt{\tilde f}_i , \vecop K(\vt{f}_j) \rangle$,
$[\mat{T_{BC}}]_{ij} \coloneqq \langle \uv n \times \vt{\tilde{f}}_i , \vecop T(\vt{\tilde{f}}_j) \rangle$,
$[\vec e_3]_i \coloneqq \langle \uv n \times \vt{f}_i, \uv n \times \vt e^\mathrm{inc} \rangle$, $[\vec h_3]_i \coloneqq \langle \uv n \times \vt{\tilde{f}}_i, \uv n \times \vt h^\mathrm{inc} \rangle$,
and $\langle \vt a, \vt b\rangle\coloneqq\int_{\Gamma_3} \vt a \cdot \vt b \diff s$.

\begin{figure}
        \centering
        \input{Figures/rwg.tex}
        \caption{Illustration of the conventions used in defining the RWG functions $\{\vt f_i\}$: the support of $\vt f_i$ is the union of the triangles $c_i^+$ and $c_i^-$, that have the RWG-defining, oriented edge $\vt e_i$ in common and are complemented by $\vr_i^+$ and $\vr_i^-$, respectively.}
        \label{fig:RWGconv}
\end{figure}
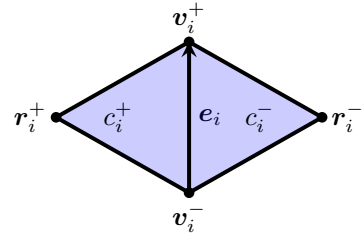

To allow for complete analyses of the matrices of the 3-dimensional problem, we introduce the loop-to-RWG and star-to-RWG transformation matrices $\mat \Lambda \in \mathbb{R}^{N \times N_l}$ and $\mat \Sigma \in \mathbb{R}^{N \times N_s}$ that are mappings between the solenoidal and RWG subspaces, and between the non-solenoidal and RWG subspaces, respectively. Here $N_s$ is the number of triangles in the discretization of $\Gamma_3$ and $N_l$ the number of vertices. These mappings can be defined using the notation illustrated in Fig.~\ref{fig:RWGconv} as
\begin{equation}
    [\mat \Sigma]_{ij} \coloneqq
    \begin{dcases}
        \phantom{-}1 & \text{if triangle $j$ is $c_i^+$,}\\
        -1 & \text{if triangle $j$ is $c_i^-$,}\\
        \phantom{-}0 & \text{otherwise,}
    \end{dcases}
\end{equation}
and
\begin{equation}
    [\mat \Lambda]_{ij} \coloneqq
    \begin{dcases}
        \phantom{-}1 & \text{if vertex $j$ is $v_i^+$,}\\
        -1 & \text{if vertex $j$ is $v_i^-$,}\\
        \phantom{-}0 & \text{otherwise.}
    \end{dcases}
\end{equation}
In addition to these standard definitions, we also introduce their counterparts acting on orthonormalized bases $\tmat \Sigma = \matRR^{-1/2} \mat \Sigma \matPP^{1/2}$ and $\tmat \Lambda = \matRR^{1/2} \mat \Lambda \matPyrPyr^{-1/2}$.
These transformation matrices can also be used to form the quasi-Helmhotlz projectors \cite{andriulli_well-conditioned_2013}
\begin{align}
    \mat P^\Sigma &\coloneqq \mat\Sigma \left(\mat\Sigma^T\mat\Sigma\right)^{+} \mat\Sigma^T\,,\\
    \matd P^\Lambda &\coloneqq \mat\Lambda \left(\mat\Lambda^T\mat\Lambda\right)^{+} \mat\Lambda^T\,,
\end{align}
and their counterparts based $\tmat \Sigma$ and $\tmat \Lambda$
\begin{align}
    \mat P^{\tilde\Sigma} &\coloneqq \tmat\Sigma \left(\tmat\Sigma^T\tmat\Sigma\right)^{+} \tmat\Sigma^T\,,\\
    \matd P^{\tilde\Lambda} &\coloneqq \tmat\Lambda \left(\tmat\Lambda^T\tmat\Lambda\right)^{+} \tmat\Lambda^T\,.
\end{align}

From the boundary element matrices $\mat Z_{2,3}$ of the normalized Calderón TE-EFIE (2D) and normalized Calderón EFIE (3D), we define the matrices $\mat C_{2,3}$
\begin{equation}
    \mat C_{2,3} \coloneqq \mat Z_{2,3} - \mat I/4\,,
\end{equation}
where $\mat I$ is the identity matrix. In the following, to simplify the notation, we will omit the subscripts $_2$ and $_3$ when the developments apply to both 2-dimensional and 3-dimensional problems.

\section{Laplacian Filters}

In this section we leverage on operator filtering results presented in \cite{merlini_laplacianFilters_2022}. 
The following notation will be used for the singular value decomposition (SVD) \cite{golub_matrix_2012} of a symmetric matrix $\mat X \in \mathbb{R}^{N_{x}\times N_{x}}$
\begin{equation}
    \mat X=\mat V_{\mat X} \mat S_{\mat X} \mat V_{\mat X}^T
\end{equation}
where $\sigma_{\mat X,i}=\left[\mat S_{\mat X} \right]_{ii}$ are the singular values  of $\mat X$.
We define the diagonal matrix $\mat S_{\mat X,n}$ as 
\begin{equation} 
  \left[\mat S_{\mat X,n} \right]_{ii} =
  \begin{dcases}
    \sigma_{\mat X,i} &  \text{if } i >  N_{x}-n\,,\\
   \phantom{-}0 & \text{otherwise,}
  \end{dcases}
  \label{eq:def_stars}
\end{equation}
and the filtered matrix $\mat X_n$ as 
\begin{equation}\label{eq:svdX}
\mat X_{n} \coloneqq \mat V_{\mat X} \mat S_{\mat X,n} \mat V_{\mat X}^T\,.
\end{equation}

\subsection{Two-Dimensional Case}
Consider the discretized two dimensional variational Laplacian $\mat L$ where
\begin{equation}
    [\mat L]_{ij}=\langle\nabla\varphi_i,\nabla\varphi_j\rangle\,,
\end{equation}
and its orthonormalized counterpart $\tmat L \coloneqq \mathh^{-1/2} \mat L \mathh^{-1/2}$.
Consistently with the general definition \eqref{eq:svdX}, we obtain the filtered Laplacian
\begin{equation}
\mat L_{n} \coloneqq \mat V_{\mat L} \mat S_{\mat L,n} \mat V_{\mat L}^T
\end{equation}
from which we can define the associated Laplacian filter
\begin{equation}\label{eq:LapFilter}
\mat P^L_{n} \coloneqq \mat L_{n}^+\mat L_{n}=\mat V_{\mat L} \mat S_{\mat L,n}^{+}\mat S_{\mat L,n} \mat V_{\mat L}^T\,.
\end{equation}
The reader should note that the SVD is used in the definitions and  will never be necessary in practice, as the FFT will be used to keep the quasi-linear complexity.
In the definitions above, the Laplacian matrix $\mat L$ should be replaced by $\tmat L$ when dealing with non-homogeneously discretized structures as is assumed in the following.

\subsection{Three-Dimensional Case}

Following \cite{merlini_laplacianFilters_2022}, we define the primal filters as
\begin{align}
    \mat P_n^\Sigma&=\mat\Sigma \left(\mat\Sigma^T\mat\Sigma\right)^{+}_n \mat\Sigma^T\,,\\
    \mat P_n^{\Lambda H}&=\mat\Lambda \left(\mat\Lambda^T\mat\Lambda\right)^{+}_n \mat\Lambda^T +\mat I-\mat P^\Sigma-\matd P^\Lambda\,,\label{eq:defLamH}
\end{align}
and the dual ones as
\begin{align}
\matd P_n^\Lambda&=\mat\Lambda \left(\mat\Lambda^T\mat\Lambda\right)^{+}_n \mat\Lambda^T\,, \\
    \matd P_n^{\Sigma H}&=\mat\Sigma \left(\mat\Sigma^T\mat\Sigma\right)^{+}_n \mat\Sigma^T+\mat I-\matd P^\Lambda-\mat P^\Sigma\,. \label{eq:defSigH}
\end{align}
The above loop and star matrices $\mat \Lambda$ and $\mat \Sigma$ should be replaced by $\tilde{\mat \Lambda}$ and $\tilde{\mat \Sigma}$
when dealing with non-homogeneously discretized structures as is assumed in the following.

\section{Filters Based Fast Direct Solver}

The filters defined in the previous section can be used to obtain a fast direct solver, in the low-frequency regime, for the two-dimensional and three-dimensional cases.
\subsection{Two-dimensional case}
Considering first the TE-EFIE formulation
\begin{equation}
     \mat Z_2 \vec j_2 = \vec v_{e,2}\,,
\end{equation}
the compact part of which is represented by the matrix $\mat C_2$, the filtered equation we propose reads
\begin{equation}
\left(\mat I/4 + \mat P^{L}_n \mat C_2 \right)\vec j_2 = \vec v_{e,2}\,.
\end{equation}
Dually, for the TE-MFIE, the compact part is just the operator $\op D$, so that the filtered TM-MFIE reads
\begin{equation}
\left(\mat I /2+\mat P^{L}_n \mat D \right)\vec j_2 = \vec v_{h,2}\,.
\end{equation}

By combining the two previous formulations, after localizing the preconditioner 
as is standard in Calderón-preconditioning and grouping the compact filtered parts, we obtain the filtered Yukawa-Calderón TE-CFIE
\begin{equation}
\left(\frac{(1+2\alpha)}{4}\mat I+ \mat P^{L}_n \left(\mat C_2+\alpha\mat D\right)\right)\vec j_2
=      \vec v_{e,2} + \alpha\vec v_{h,2}\,.
\end{equation}

\subsection{Three-dimensional case}

Similarly, in the 3D case, the filtered EFIE we propose is
\begin{equation}
\left(\mat I/4 + (\mat P_n^{\Lambda H} + \mat P_n^{\Sigma}) \mat C_3\right) \vec j_3 = \mat v_{e,3}\,,
\end{equation}
and, because its compact part is just the operator $\mat{\tilde{K}}$, the filtered MFIE reads
\begin{equation}
\left(\mat I/2 + (\mat P_n^{\Lambda H} + \mat P_n^{\Sigma}) \mat{\tilde{K}} \right) \vec j_3 =  \vec v_{h,3}\,.
\end{equation}

By combining the previous two formulations, after localizing the preconditioner 
as is standard in Calderón preconditioning, and grouping the compact filtered parts, we obtain the filtered Yukawa-Calderón CFIE
\begin{multline}
\left(\frac{(1+2\alpha)}{4}\mat I + (\mat P_n^{\Lambda H} + \mat P_n^{\Sigma}) \left( \mat C_3 + \alpha \mat{\tilde{K}} \right) \right) \vec j_3\\ = \vec v_{e,3} + \alpha \vec v_{h,3}\,.
\end{multline}

\subsection{Inversion of the structured matrices}

Each of the above equations is in the form $\beta \mat I + \mat C$, where $\beta$ is a constant and $\mat C$ is a low rank matrix that can be expressed in the form $\mat C\eqqcolon\mat U\mat V^T$ by using a skeleton-detecting algorithm \cite{bebendorf_hierarchical_2008,liberty_randomized_2007}. The peculiar single-skeleton structure we obtain in this way allows for a linear-in-complexity direct inverse in the form 
\begin{multline}
    \left(\beta\mat I + \mat C\right)^{-1}=\left(\beta \mat I + \mat U\mat V^T\right)^{-1}=\\
    \mat I \beta^{-1} - \beta^{-2}\mat U \left( \mat I_r + \beta^{-1}\mat V^T \mat U \right)^{-1} \mat V^T\,,
\end{multline}
where $\mat I_r$ is the identity matrix that has as dimension the rank of $\mat C$, using the Woodbury matrix identity \cite{henderson_deriving_1981}.

\begin{table*}[t]
\caption{Compression efficiency and solution accuracy obtained with the perturbed circle with a compression tolerance of $\epsilon = \num{e-3}$ and filtering point of $200$.}
\label{table}
\centering
\begin{tabular}{cccccc}
\hline
Number of unknowns & 1004 & 2008 & 4016 & 8032 \\
\hline
Solution relative error & $\num{0.23e-3}$& $\num{0.24e-3}$& $\num{0.24e-3}$& $\num{0.24e-3}$\\
\hline
Memory usage without compression (MBytes) & \num{16}& \num{64}& \num{258}&  \num{1032}\\
\hline
Memory usage skeleton (MBytes) & \num{2.37}& \num{1.54}& \num{3.08}& \num{6.16}\\
\hline
\end{tabular}
\end{table*}

\section{Numerical Results}

\begin{figure}
\centering
\includegraphics[width=1\columnwidth]{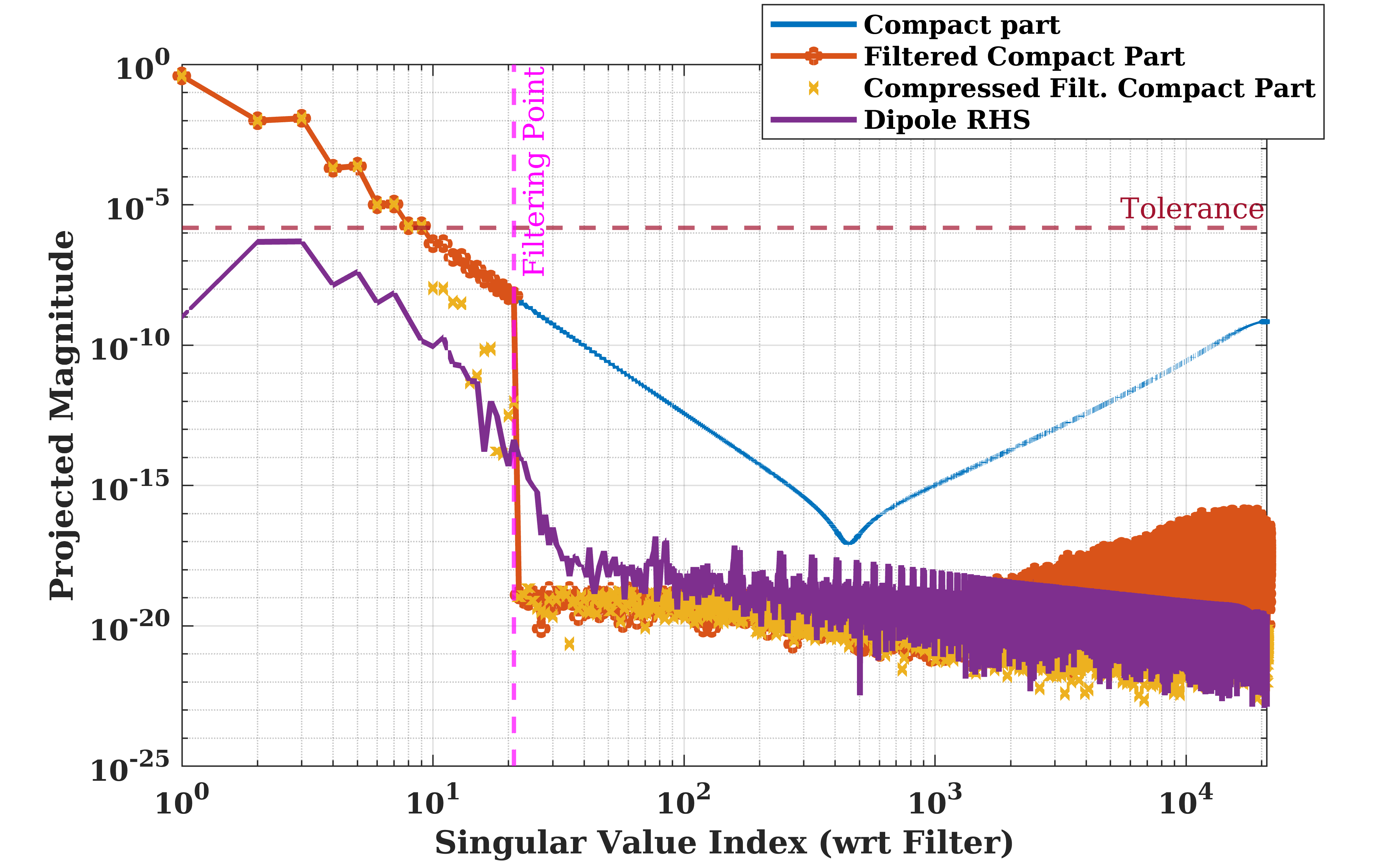}
\caption{Projections of the compact part ($\mat C$), the filtered compact part ($\mat{\tilde{C}}$), and the dipole RHS on the singular vectors of the filter $\mat P^L_N$. The yellow markers indicate which singular values/singular vectors are kept after compression.}
\label{fig:spectra}
\end{figure}

We first evaluated the fast direct solution scheme on an ellipse characterized by a semi-major axis of \SI{1.42}{\meter} and a semi-minor axis of \SI{1.32}{\meter} and illuminated by a dipole source radiating at $k=\SI{0.4}{\radian \per \meter}$. Figure~\ref{fig:spectra} shows the projections of the compact part ($\mat C$), the filtered compact part ($\tmat{C} = \mat P^L_n \mat C$), the compressed and filtered compact part ($\mat U \mat V^T$), and the dipole right hand side (RHS) $\vec v$ on the singular vectors of the filter $\mat P^L_N$. The projection of $\mat C$ on this basis makes apparent the spectral deviation induced by the discretization error (i.e. growth of the projected magnitude for increasing singular value index) which prevents its proper compression using standard techniques. The projection of $\mat{\tilde{C}}$ is free from this shortcoming, thanks to the filtering which allows for the compression of $\mat{\tilde{C}}$ at any tolerance above the last singular value filtered out (see Fig.~\ref{fig:spectra}). Moreover, the projected RHS in Fig.~\ref{fig:spectra} is band-limited with respect to the filter, which ensures that the filtering does not deteriorate the precision of the solution.

\begin{figure}
\centering
\includegraphics[width=1\columnwidth]{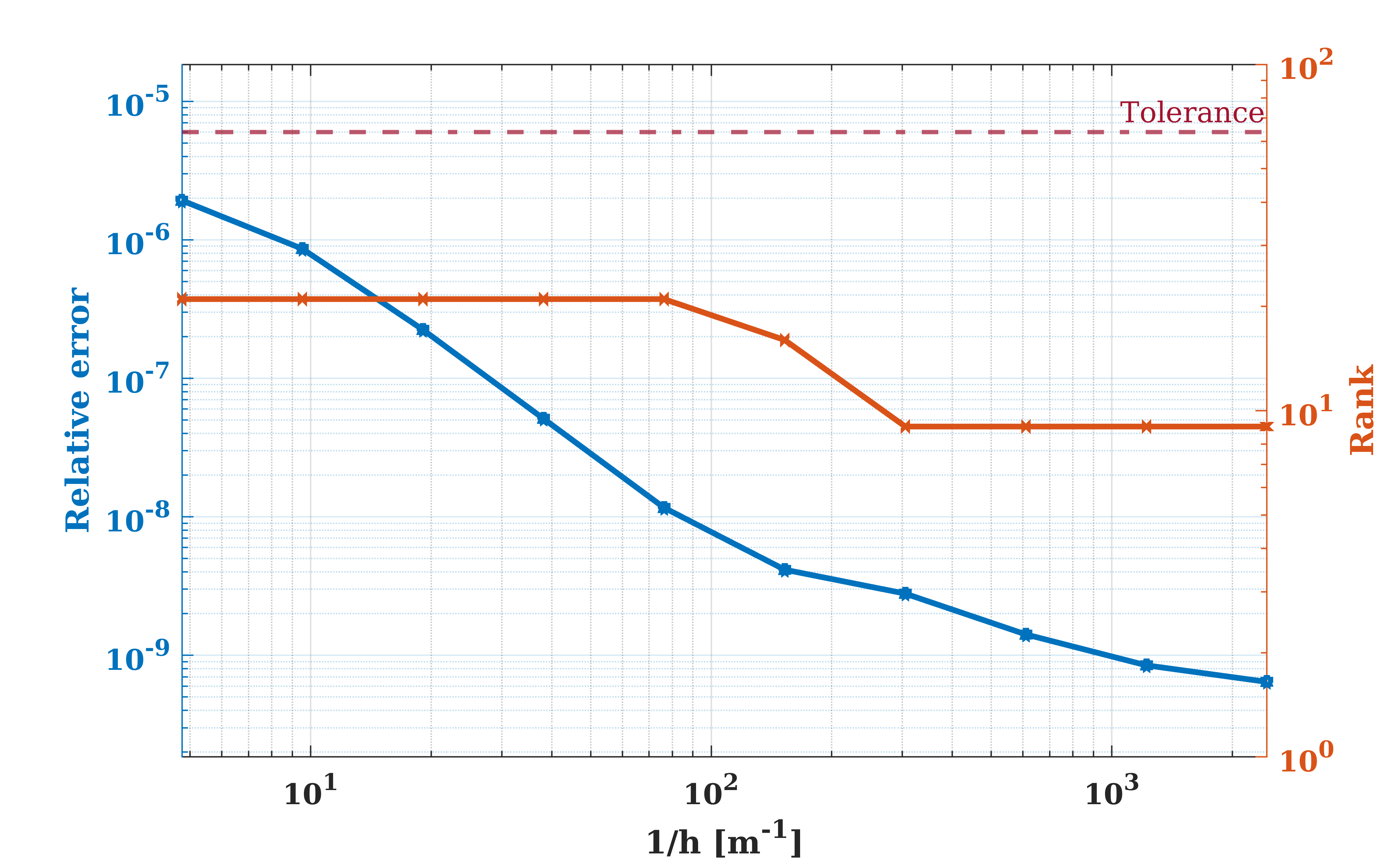}
\caption{Accuracy of the solution and inner rank of the skeleton as a function of the inverse edge length, with a compression tolerance of $\epsilon = \num{6e-6}$ and a filtering point of $21$.}
\label{fig:rank_re_disc}
\end{figure}

To demonstrate the accuracy of the solution and the compression efficiency, we show in Fig.~\ref{fig:rank_re_disc} the error on the solution and the inner rank of the skeleton as a function of the inverse edge length of the discretized ellipse. For increasing refinement, the error with respect to the non compressed operator decreases and the rank saturates to a constant, which validates the effectiveness of the scheme.

Finally, the ellipse has been substituted by a perturbed circle with radius $r(\theta)=2 + 0.2 \sin(8\theta)$. As shown in Table~\ref{table}, the storage requirement for the skeleton, from a given refinement, grows linearly with the number of unknowns and, at the same time, the error on the solution remains below the desired tolerance. 

\section*{Acknowledgment}

This work has been funded in part by the European Research Council (ERC) under the European Union’s Horizon 2020 research and innovation program (ERC project 321, grant No.724846), in part by the Italian Ministry of University and Research within the Program PRIN2017, EMVISIONING, under Grant 2017HZJXSZ, CUP:E64I190025300, and in part by the ANR Labex CominLabs under the project ``CYCLE''.


\bibliographystyle{IEEEtran}
\bibliography{MyLibrary,bib_adrien}




%

\end{document}

%% file: Macros.tex
\newcommand{\vt}[1]{\boldsymbol{#1}} 
\newcommand{\vr}{\boldsymbol{r}} 

\newcommand{\diff}{\,\mathrm{d}} 

\newcommand{\veg}[1]{\bm{#1}}               
\newcommand{\vet}[1]{\bm{\mathfrak{#1}}}    
\newcommand{\mat}[1]{\mathbf{#1}}       
\newcommand{\matPyrPyr}{\mathbf{G}_{\lambda \lambda}}
\newcommand{\matPP}{\mat G_{pp}}
\newcommand{\matRR}{\mat G_{f f}}
\newcommand{\matMix}{\tmat G_{f \tilde{f}}}
\newcommand{\matMixTr}{\tmat G_{\tilde{f} f}}
\newcommand{\mathh}{\mat G_{\varphi \varphi}}

\newcommand{\matd}[1]{\mathbb{#1}}   
\newcommand{\uv}[1]{\hat{\veg{#1}}}         
\renewcommand{\vec}[1]{\mathsfbfit{#1}}     
\newcommand{\op}[1]{\mathcal{#1}}           
\newcommand{\vecop}[1]{\bm{\mathcal{#1}}}   

\newcommand{\im}{\mathrm{i}}  

\newcommand{\R}{\mathbb{R}}

\newcommand{\tmat}[1]{\mat{\tilde{#1}}}

%% file: Packages.tex
\usepackage[utf8]{inputenc}

\usepackage[T1]{fontenc}
\usepackage{float}
\usepackage{tikz}
\usetikzlibrary{arrows,matrix,shadows,positioning}
\usepackage{pgfplots}
\usepackage{pgfplotstable}
\usepackage{siunitx}
\usepackage{cleveref}
\pgfplotsset{compat=newest}
\usepackage{multirow}

\usepackage{amsthm}
\DeclareMathAlphabet{\mathbfsf}{\encodingdefault}{\sfdefault}{bx}{n}
\usepackage{bm}
\usepackage[OMLmathsfit]{isomath}


%% file: Figures/rwg.tex
\definecolor{linec}{rgb}{0,0,0}
\definecolor{trianglec}{rgb}{0,0,255}
\definecolor{nodec}{rgb}{0,0,0}

\begin{tikzpicture}[line cap=round,line join=round,>=triangle 45]
\coordinate[label=right:$\vt e_i$] (center) at (0,0);
\coordinate[label=below:$\vt v_i^-$] (vm) at (0,-1);
\coordinate[label=above:$\vt v_i^+$] (vp) at (0,1);
\coordinate[label=right:$\vt r_i^-$] (rm) at (1.75,0);
\coordinate[label=left:$\vt r_i^+$] (rp) at (-1.75,0);

\draw [line width=1.5pt,color=linec,fill=trianglec,fill opacity=0.2] (vm) -- (rp) -- (vp) -- (rm) -- (vm);
\draw [line width=1.5pt,color=linec,-stealth] (vm) -- (vp);

\fill [color=nodec] (vm) circle (2pt);
\fill [color=nodec] (vp) circle (2pt);
\fill [color=nodec] (rm) circle (2pt);
\fill [color=nodec] (rp) circle (2pt);

\coordinate[label=left:$c_i^-$] (tm) at (1.25,0);
\coordinate[label=right:$c_i^+$] (tp) at (-1.25,0);
\end{tikzpicture}